\documentclass[a4paper,12pt]{amsart}

\usepackage[margin=3cm, headsep=0.8cm]{geometry}
\usepackage[usenames,dvipsnames]{xcolor}

\usepackage{graphicx}

\usepackage{tikz-cd}

\newtheorem{thm}{Theorem}

\newtheorem{conj}{Conjecture}

\newtheorem{claim}[thm]{Claim}
\newtheorem{lem}[thm]{Lemma}
\newtheorem{slem}[thm]{Sublemma}
\newtheorem{cor}[thm]{Corollary}
\newtheorem{prop}[thm]{Proposition}

\newtheorem{definition}[thm]{Definition}

\newcommand{\N}{\mathbb N}
\newcommand{\R}{\mathbb R}

\def\al{\alpha}
\def\ga{\gamma}
\def\Ga{\Gamma}

\def\eps{\epsilon}

\def\Si{\Sigma}
\def\om{\omega}

\def\3{\ss}

\def\wlim{\mathop{\hbox{$\om$-lim}}}

\def\D{\partial}

\def\spt{\operatorname{spt}}

\def\<{\langle}
\def\>{\rangle}

\newcommand{\bcl}{\begin{claim}}
\newcommand{\ecl}{\end{claim}}
\newcommand{\bcor}{\begin{cor}}
\newcommand{\ecor}{\end{cor}}

\newcommand{\bdfn}{\begin{definition}}
\newcommand{\ben}{\begin{enumerate}}
\newcommand{\bit}{\begin{itemize}}
\newcommand{\blem}{\begin{lem}}
\newcommand{\bslem}{\begin{slem}}
\newcommand{\bprop}{\begin{prop}}
\newcommand{\bthm}{\begin{thm}}
\newcommand{\edfn}{\end{definition}}
\newcommand{\een}{\end{enumerate}}
\newcommand{\eit}{\end{itemize}}
\newcommand{\elem}{\end{lem}}
\newcommand{\eslem}{\end{slem}}
\newcommand{\eprop}{\end{prop}}
\newcommand{\ethm}{\end{thm}}

\title[Curvature bounds in 2D]{Curvature bounds of subsets in dimension two}
\author{Alexander Lytchak, Stephan Stadler}

\newcommand{\Addresses}{{\bigskip\footnotesize
\noindent Alexander Lytchak, 
\par\nopagebreak\noindent\textsc{Mathematisches Institut, Universit\"at K\"oln, Weyertal 86 - 90, 50931 K\"oln, Germany}
\par\nopagebreak
\noindent\textit{Email}: \texttt{alytchak@math.uni-koeln.de}

\medskip

\noindent Stephan Stadler,
\par\nopagebreak\noindent\textsc{Max Planck Institute for Mathematics, Vivatsgasse 7, 53111 Bonn, Germany}
\par\nopagebreak
\noindent\textit{Email}: \texttt{stadler@mpim-bonn.mpg.de}

}}



\begin{document}

\begin{abstract}
We show that closed subsets with vanishing first homology in two-dimensional spaces
inherit the upper curvature bound from their ambient spaces and discuss topological applications.
\end{abstract}

\maketitle

\section{Introduction}

This paper concerns the intrinsic geometry of subsets in two-dimensional metric spaces with upper curvature bounds.
The main geometric result is

\bthm\label{thm_intcurv_intro}
Let $X$ be a two-dimensional  contractible CAT($\kappa$) space.
 Let $A\subset X$ be a closed, Lipschitz connected  subset with $H_1(A)=0$.
 Then $A$ is  a CAT($\kappa$) space with respect to the induced intrinsic metric.
\ethm

For $\kappa =0$,  this confirms a  folklore
conjecture, which appeared in print in \cite[Conjecture~1]{Ber}.
	Related statements and conjectures about subsets of non-positively curved spaces can be found in \cite[Chapter~4]{AKP_inv}.

	Some special cases  of  Theorem~\ref{thm_intcurv_intro} are  known. In \cite{Bishop} and later in  \cite{Ricks_planar},
	it was shown that Jordan domains in the euclidean  plane are
	CAT(0).	A more general version appeared in \cite{LWcurv}.
	In \cite{Ricks}, Theorem~\ref{thm_intcurv_intro}  is proved  for CAT(0)
	euclidean simplical complexes.
	Another special case plays a central role in \cite{NSY}.

The contractibility assumption is redundant for $\kappa \leq 0$; for $\kappa >0$ it is satisfied if the diameter of $X$ is less than $\frac \pi {\sqrt \kappa}$.  For  $\kappa>0$ the statement is wrong without the contractibility assumption:  A closed metric ball
	of radius $\pi >r>\frac \pi 2$ in the round sphere $\mathbb S^2$ is contractible but not CAT(1) in its intrinsic metric.

Localizing the above result  we deduce the following:
\begin{cor}\label{cor_local}
	Let $Y$ be a
	 metric
	 space of curvature bounded above by $\kappa$ and dimension two. Let $A\subset Y$ be  closed,   Lipschitz connected and  locally simply connected.
Then $A$ has curvature bounded above by $\kappa$ with respect to the intrinsic metric.
\end{cor}

If $\kappa \leq 0$,
 the theorem of Cartan--Hadamard
	 implies that
	the universal covering of $A$ is contractible.
Hence $A$ is \emph{aspherical} in the sense that all higher homotopy groups of $A$ vanish:

\bcor\label{cor_asph_intro}
Let $Y$ be a two-dimensional space of non-positive curvature and let $A\subset Y$ be a closed,  Lipschitz connected and locally simply connected  subset.
Then $A$ is aspherical with respect to the topology induced by the intrinsic metric.
\ecor

Somewhat surprisingly, no topological assumption is needed for our next conclusion. Indeed, we
obtain the following topological statement about \emph{all} subsets of non-positively curved
metric spaces of dimension two.

\bthm\label{thm_lipasph_intro}
Let $Y$ be non-positively curved and  two-dimensional.
Let $A\subset Y$ be   an arbitrary subset.
Then all higher Lipschitz homotopy groups of $A$ vanish:
Every Lipschitz $n$-sphere in $A$ with $n\geq 2$  bounds a Lipschitz ball in $A$.
\ethm

This result is of geometric origin and is deduced from our main theorem.
It has the following purely topological application:

\bcor\label{cor_opasph_intro}
Let $Y$ be a  two-dimensional space of non-positive curvature. Then any neighborhood retract $A\subset Y$
is aspherical.
\ecor

It is known, but surprisingly difficult to prove  that \emph{all} subsets of the euclidean plane are aspherical, \cite{CCZ}. This and the results above make the following generalization of the famous Whitehead Conjecture \cite{Whitehead}  plausible:

\begin{conj}
 Let  $X$ be a two-dimensional aspherical space.  Then any subset $A$ of $X$ is aspherical.
\end{conj}

 Possibly, a combination of the geometric ideas of the present paper and the purely topological methods of \cite{CCZ} may lead to the  resolution of this conjecture
for non-positively curved spaces.

Since there is no assumption on local compactness in Theorem~\ref{thm_intcurv_intro}, we gain some information on coarse topology
of high-dimensional spaces.
Recall that a CAT(0) space has  {\em asymptotic rank at most two},  if no asymptotic cone contains an isometric copy of euclidean $3$-space
\cite{G_asym, CaLy,  W_asym}.

\bprop\label{prop_coarsetop}
Let $X$ be a CAT(0) space of asymptotic rank at most two and let $A\subset X$ be arbitrary.
Then for every $n\in\N$ and $\eps>0$ there exists $L_0>0$ such that the following holds for all $L\geq L_0$.
Any $L$-Lipschitz sphere $f:\mathbb S^n\to A$ bounds a Lipschitz ball in $N_{\eps L}(A)$.
\eprop

As the proof shows, the numbers $\epsilon, L_0$ are independent of $A$, but only depend on $X$. Moreover, the map $f$ extends to a  $(\pi L)$-Lipschitz map $F$ from
	the euclidean unit  ball $B^{n+1}$  into the $(\epsilon L)$-neighborhood of the image of $f$.  

Proposition \ref{prop_coarsetop}  applies, in particular, to arbitrary subgroups of rank two CAT(0) groups, finitely presented or not, compare the discussion on the Coarse Whitehead Conjecture in \cite[p.~26]{Ka_problems}.

We expect our results to simplify the description of geodesically complete two-dimensional CAT($\kappa$) spaces obtained and announced in  \cite{NSY}.   Moreover, we expect the results to  facilitate a good understanding of two-dimensional CAT($\kappa$) spaces beyond geodesic completeness.  For instance, they  might lead to a resolution of  the following conjecture of potential relevance to geometric group theory:

 \begin{conj}
  Any compact two-dimensional non-positively curved space is homotopy equivalent to a finite, two-dimensional, non-positively curved euclidean complex.
 \end{conj}

	We want to point out that all the results above trivially hold in dimension one, since any one-dimensional CAT($\kappa$) space is covered by a tree. On the other hand, all results completely  fail   in dimension at least three:  already the  complement of an open ball in $\R^3$ is not non-positively curved and not aspherical.

\bigskip

\noindent{\bf Sketch of proof.}
In order to control the curvature bound of $A$ we need to majorize arbitrary  Jordan curves $\Gamma \subset A$  by some CAT($\kappa$)-discs inside $A$. We use the homological assumption, to find
a geodesic of $X$ completely contained in $A$ which subdivides $\Gamma$ into two smaller Jordan curves. Iterating this process, we  subdivide $\Gamma$ into a collection of $2^k$ smaller Jordan curves.
 The technically most challenging part of the proof controls this cutting process and confirms that the arising new Jordan curves  will have arbitrary small diameters after sufficiently many steps. Now we majorize the small Jordan curves inside  the ambient space $X$ and observe that these majorizations glue together to a majorization of $\Gamma$ within  a   small neighborhood of $A$ in $X$.  
A limiting argument provides the required majorization contained in $A$.

\bigskip

\noindent{\bf Acknowledgments.}  We are grateful to Anton Petrunin for helpful comments.
Both authors were  supported by DFG grant SPP 2026.

\section{Metric geometry}

\subsection{Basics and notation}

We refer
to \cite{BBI}, \cite{Ballmann}, \cite{AKP} for background on metric geometry and
CAT($\kappa$) spaces. Let us summarize notation and basic facts.
As usual $\R^n$  will denote the euclidean space and  $\mathbb S^{n-1}\subset\R^n$ the unit sphere.

The distance  on a metric space  $X$ will be denoted by $|\cdot,\cdot|_X$ and if there is no risk for
confusion by $|\cdot,\cdot|$. If $A\subset X$ is a subset, then we denote  by $N_r(A)$ and $\bar N_r(A)$ its open respectively closed $r$-neighborhoods. If $A$ is just a point $x$, then $N_r(A)$  is  the open
$r$-ball which we
denote by
$B_r(x)$. The closed $r$-ball will be denoted by $\bar B_r(x)$.

The length of a curve $c$ in a metric space  $X$ will be denoted  by $\ell (c) \in [0, \infty]$. The space $X$ is called Lipschitz connected if any two points in $X$ are joined by a curve of finite length.
$X$ is called an
{\em intrinsic space},
 if the distance between any two points
is equal to the greatest lower bound for lengths of curves connecting those points.

Isometric embeddings of intervals will be called {\em geodesics} and in the case of compact intervals also {\em geodesic segments}
or simply {\em segments}. If $c$ is a geodesic segment, then we will denote its boundary points by $\D c$.
 The space $X$ itself will be called {\em geodesic} if any two points in $X$
are joined by a geodesic.

Any Lipschitz connected metric space has a canonical \emph{induced intrinsic metric},  \cite[2.3.3]{BBI}.  The length of all curves for the given metric and for the induced intrinsic metric coincide.

A {\em triangle} in $X$ consists of three points and three geodesics connecting them.
The three geodesics  are called the {\em sides} of the triangle.
 For $\kappa\in\R$, let $D_\kappa\in(0,\infty]$ be the diameter of the complete, simply connected model surface $M^2_\kappa$
of constant curvature $\kappa$. For any triangle $\triangle$ with perimeter $<2 D_\kappa$ in a metric space $X$, we can find a \emph{comparison triangle} $\tilde\triangle\subset M^2_\kappa$
such that corresponding sides have equal lengths.

 A complete metric space $X$ is CAT($\kappa$) if points at distance $<D_\kappa$ in $X$
are joined by geodesics and if all triangles in $X$ with perimeter $<2 D_\kappa$ are \emph{not thicker} than their comparison triangles.  A metric space  $X$ is said to have {\em curvature bounded above} by $\kappa$, if $X$ is  locally
CAT($\kappa$).

 Any CAT(0) space is contractible; the universal covering of any complete non-positively curved space is CAT(0),   \cite[8.13.1, Theorem of Cartan--Hadamard]{AKP}.

In a CAT($\kappa$) space $X$ any point $y$ in $B_{D_{\kappa}} (x)$ is connected with  $x$ by a unique geodesic, denoted by $xy$ and  depending continuously on $y$.
 Angles between geodesics starting at  $x$ are well defined.
There is a {\em space of directions} $\Si_x X$ which is CAT(1) with respect to the angle metric. The logarithm map
\[\log_x:B_{D_\kappa}(x)\setminus\{x\}\to \Si_x X \,,\]
which assigns to  $y$ the starting direction of $xy$
 is a homotopy equivalence \cite{Kramer}.

In  a CAT($\kappa$) space $Z$, there are natural ways of straightening singular simplices using iterated geodesic coning or barycentric simplices, see \cite[Section~6.1]{KleinerLeeb} or \cite[Lemma 5.1]{Kleiner}.
As  a consequence,  for any finite simplicial complex $K$, Lipschitz maps are dense in the space of continuous maps $K\to Z$ with respect to compact-open topology.
Moreover, if a continuous map $f$ is already Lipschitz continuous on a subcomplex $K_0$ of $K$ then the Lipschitz continuous approximations $\tilde f$ of $f$  can be chosen to coincide with $f$ on $K_0$.

On a space $X$ of curvature bounded above, there is a natural notion of dimension $\dim (X)$,
  introduced and investigated by
Kleiner in \cite{Kleiner}.
It is equal to the supremum of topological dimensions of compact subsets of  $X$, satisfies $\dim (X)= \sup _{x\in X} \{ \dim (\Sigma _xX) +1 \}$ and coincides with the homological dimension of $X$.

The homological dimension of a space $X$ is the supremum of all $n$, such that
	for some open pair $V\subset U \subset X$ the relative  homology  $H_n(U,V)$ is non-zero.  Here and below  $H_n$ denotes singular homology with $\mathbb Z$ coefficients.

\section{Majorizations}

A \emph{Jordan curve} in a metric space $X$ is a subset
homeomorphic to a circle.  We say that a metric space $Y$ majorizes a  rectifiable Jordan curve $\Gamma$ in a metric space $X$
if there exists a $1$-Lipschitz map $f:Y\to X$ which sends a Jordan curve $\Gamma'  \subset Y$ bijectively in an arc length preserving way onto $\Gamma$.
By Reshetnyak's Majorization Theorem, \cite[8.12.4]{AKP}, any  Jordan curve of length $< 2D_{\kappa}$ in a CAT($\kappa$) space is majorized by a closed convex subset of $M^2_{\kappa}$. On the other hand, we have

\begin{prop}\cite[Proposition~4.2]{LS_improv}  \label{prop_major}
Let $X$ be a complete intrinsic space.  If any Jordan curve $\Gamma$  of length $<2D_{\kappa}$ in $X$
is majorized by some CAT($\kappa$) space $Y_{\Gamma}$, then $X$ is CAT($\kappa$).
\end{prop}

Majorizations stay close to the curve:

\blem\label{lem_controlledmaj}
Let $X$ be a CAT($\kappa$) space and let $\Ga\subset X$ be a Jordan curve of length $<2D_{\kappa}$.
Suppose that the diameter of $\Ga$ is at most
$\eps \leq \frac {D_{\kappa}} 2$.
Then $\Ga$ has a majorization by a convex subset of $M^2_{\kappa}$
whose  image  in $X$ is  contained in $\bar N_{\epsilon} (\Gamma)$.
\elem

\proof
For any $x\in \Gamma$, the closed ball $\bar B_{\eps} (x)$ is CAT($\kappa)$ and contains $\Gamma$.
Hence, we can find a majorization of $\Gamma$ within this ball.
\qed

\medskip

The following simple gluing/subdivision lemma will be repeatedly applied in the proof of the main theorem.

\blem\label{lem_glumaj}
Let $X$ be a CAT($\kappa$) space and $\Ga^\pm\subset X$ rectifiable Jordan curves.
Suppose $\Ga^+$ intersects $\Ga^-$ in a geodesic segment $\ga$ and denote by $\Ga$ the Jordan curve $(\Ga^+\cup\Ga^-)\setminus\ga$.
Let $Z^\pm$ be CAT($\kappa$) discs
and let $w^\pm:Z^\pm\to X$ be majorizations of $\Ga^\pm$ which restrict to arc length preserving homeomorphisms $\D Z^\pm\to\Ga^\pm$.
Then there is a CAT($\kappa$) disc $Z$ and a majorization $w:Z\to X$ of $\Ga$ which sends $\D Z$ in an arc length preserving way onto $\Ga$ and whose image
is  the union of the images of $w^\pm$.
\elem

\proof
Since $\ga$ is a geodesic in $X$ and $w^\pm$ is $1$-Lipschitz, the subarc $c^\pm$ of $\D Z^\pm$ which gets mapped to $\ga$ has to be geodesic in $Z^\pm$,
cf. \cite[8.12.2]{AKP}.
Let $f:c^+\to c^-$ be the canonical isometry. By Reshetnyak's gluing theorem~\cite[8.9.1]{AKP}, the space $Z:=Z^+\cup_f Z^-$ is  CAT($\kappa$). Moreover, $Z$ is homeomorphic to a closed disc and contains $Z^\pm$ as convex subspaces. We define $w:Z\to X$ such that it restricts to $w^\pm$ on $Z^\pm$. By construction, $w$  is a well-defined majorization of $\Gamma$, as required.
\qed

\section{Support sets  of cycles}

Support sets of top-dimensional cycles provide a useful tool in the study of finite dimensional spaces \cite{KleinerLeeb, BKS, H, Ricks}.
We will recall their definition here and prove basic properties required in later sections, in order to define our cutting procedure for
Jordan curves. Support sets were also used by Ricks in his proof of Theorem~\ref{thm_intcurv_intro} for simplicial complexes, although in a different vein.

\bdfn
Let $Y$  be a subset of a metric space $X$. The {\em support} $\spt(\al)$ of a homology class $\al\in H_n(X,Y)$ is the set of all points $x\in X\setminus Y$ such that  the image of $\alpha$ is non-trivial under  the inclusion homomorphism
\[H_n(X,Y)\to H_n(X,X\setminus  \{x\})\]
\edfn

The support of
 $\al$
is  the intersection of $X\setminus Y$ with all images of chains representing $\al$. Thus, $\spt (\al)$ is closed in $X\setminus Y$ and its closure in $X$ is compact.

For instance, the support of the fundamental class $[M] \in H_n(M)$ of a compact, oriented
   $n$-manifold  $M$  is the whole manifold $M$.

The following result has been verified in \cite[p. 2342]{H}:
\begin{lem} \label{lem_new}
	Let $Y$ be a subset of a metric space $X$ of homological dimension $n$.  Let $S$ be the support of a class $\alpha \in H_n (X,Y)$.  Then, for any neighborhood $U$ of $S\cup Y$, the class $\alpha$ can be represented in $U$; thus, $\alpha$ is in the image of $i _{\ast}:H_n(U,Y) \to H_n (X,Y)$.
\end{lem}

If $X$ is contractible, then  the boundary homomorphism $\partial :H_n(X,Y) \to H_{n-1} (Y)$ is an isomorphism and
the support of  $\alpha \in H_n(X,Y)$ is the set of all   $x\in X\setminus Y$
such that the image  $i_{\ast} (\partial \alpha) \in H_{n-1} (X\setminus \{x \})$ is non-zero  for  the inclusion $i:Y\to X\setminus \{x \}$.

\begin{cor} \label{cor_fin}
	Let $X$ be   a contractible metric space of homological dimension $n$.
	Let $Y\subset X$ and  $\alpha \in H_n (X,Y)$ be given.  Set  $S=\spt (\alpha) \subset X\setminus Y$.
	Let $p\in S$ and $0<r <  |p,Y|$  be arbitrary. Then, for  all neighborhoods $W$ of $\partial B_r (p) \cap S$
	 in $X\setminus \{p \}$,   the canonical map  $i_{\ast}: H_{n-1}(W ) \to H_{n-1} (X\setminus \{p\})$ is non-trivial.
\end{cor}

\proof
 Fix a neighborhood $W$, which we may assume to be open.  Assume that  we find  a larger neighborhood $V\supset W$ such that  the non-zero  image of $\alpha \in H_n (X,Y)$  in $H_n (X, X\setminus \{p \})$ can be represented by an element
 $\beta  \in H_n (V,W)$.     Then, using that the connecting homomorphism $\partial  :H_n(X,X\setminus \{p\} ) \to H_{n-1} (X\setminus\{p\})$ is
 an isomorphism, we would deduce that  $i_{\ast}  (\partial \beta)  $ is non-zero in $H_{n-1} (X\setminus \{p \})$.

 In order to find such $V$,  consider $\hat V:= B_r (p) \cup W \cup (X\setminus \bar B_r(p))$ and its subset
 $\hat W:= W\cup (X\setminus \bar B_r(p))$.  By Lemma~\ref{lem_new}, the class $\alpha$ can be represented in $H_n(\hat V,  Y)$.     Thus, the image of $\alpha$ in $H_n (X,X\setminus \{p \})$
 can be represented by an element  $\hat \alpha \in H_n (\hat V, \hat W)$.   Now we use excision and see that $\hat \alpha$ can be represented by an element in $H_n (V,W)$ where
 $V= \hat V\setminus (\hat W\setminus W)=B_r (p) \cup W$.  This finishes the proof.
\qed

\medskip

The following observation is essentially contained in  \cite[Lemma  2.4]{Ricks}.

\blem\label{lem_acc}
Let  $X$ be
contractible metric space of homological dimension $n$.
Let $M\subset X$ be a
  compact oriented $(n-1)$-manifold with fundamental class   $[M] \in H_{n-1} (M)$.
  Consider the unique  $\alpha  \in  H_n(X,M)$ with $\partial \alpha=[M]$ and 
  let $S\subset X\setminus M$  be the support of  $\alpha$.
 Then  $S \neq \emptyset $ and  $\overline{S}=S\cup M$.
\elem

\proof
Let $x\in M$ and $\delta >0$ be arbitrary.  We claim that there exists an open neighborhood
$O$ of $M$, such that $O\cup B_{2\delta } (x) =X$ and such that the image  $\hat \alpha$ of
$\alpha$ in $H_n (X,O)$ is non-zero.  Once the claim is verified,  the support $\hat S$ of $\hat \alpha$ would be a non-empty subset of $X\setminus O$, since otherwise  an application of Lemma~\ref{lem_new} with  $U=O$ would provide a contradiction to $\hat \alpha \neq 0$.   Any $p\in \hat S$ would also be contained in the support $S$ of $\alpha$.   Since $\delta$ was arbitrary, this would imply  $x\in \bar S$. But since $x$ was arbitrary, this would show $M\subset \bar S$.  On the other hand, $S$ is closed in $X\setminus M$, thus, we would deduce $\bar S= S\cup M$.  The statement $S\neq \emptyset$ would follow as well.

 In order to verify the claim, we use that
 $M$ is an absolute neighborhood  retract, and  find a
  retraction $r:V\to M$ of a neighborhood $V$ of $M$, \cite{Hanner}.   Restricting to  a smaller
 neighborhood, if needed, we may assume that $T:=r^{-1} (x)$ has diameter smaller $\delta$.
Then the injectivity of the map
$H_{n-1} (M)\to H_{n-1} (M,M\setminus \{x\})$ implies that the image of $[M]  \in
H_{n-1} ( V)$ is  non-zero
in
$H_{n-1} ( V,  V\setminus T)$.

  Setting $O:= V\cup (X\setminus \bar B_{\delta}  (x)))$ we deduce by excision
 \[H_{n-1} (O, O\setminus T) = H_{n-1} ( V,  V\setminus T) \;.\]
  Thus,  $[M]$ is non-zero in   $H_{n-1} (O, O\setminus  T)$, hence also in $H_{n-1} (O)$.  Therefore, $[M]$ is not in the image of  $\partial :H_n (O,M)\to H_{n-1} (M)$.  Now, the long exact sequence of the triple $(M,O,X)$ shows that $\hat \alpha \neq 0$.
\qed

\medskip

If $X$ is a CAT($\kappa$) space of diameter $<D_{\kappa}$, the logarithm map $\log _x :X\setminus \{x \} \to \Sigma _x X$  is a homotopy equivalence. Thus, for a subset $Y\subset X$ and a class $\alpha \in H_n (X,Y)$,  a point $x \in X\setminus Y$ is in the support of $\alpha$  if and only if $(\log_{x}) _{\ast} (\partial \alpha) \neq 0 \in H_{n-1} (\Sigma _x X)$,
compare \cite[Definition~1.3]{Ricks}, 
\cite[p.2345]{H}.

\begin{prop}\cite[Lemma~3.1]{BKS}\cite[Lemma~A-11]{H}\label{prop_geoext}
Let $X$ be an $n$-dimen-sional CAT($\kappa$) space of diameter $<D_{\kappa}$, let   $Y\subset X$  be a closed subset and $\al$ a class in $H_n(X,Y)$ with support $S:=\spt(\al)$.
Then the following {\em geodesic extension property} holds. For any pair of points $x\in X$ and $p\in S$
the segment  $xp$ extends beyond $p$ to a point $y\in Y$ such that the subsegment $py$  lies in $S\cup\{y\}$.
\end{prop}

\proof

It is enough to find, for any  $p\in S$, $x\in X\setminus \{ p \}$  and any $0 < r<  |p,Y|$, a point $q\in \partial B_r (p) \cap S$
such that $p$ lies on the segment $xq$.  Then an iteration of this property,
as at the end  of  \cite[Lemma~A-11]{H}, finishes the proof.

If there is no such  $q$, then there exists a neighborhood $W$ of  $\partial B_r (p) \cap S$ in $X$, such that
$p$ does not lie on geodesic from $x$ to a point in $W$.  Thus, the geodesics towards $x$ provide a contraction of $W$ inside $X\setminus \{ p \}$, in contradiction to Corollary \ref{cor_fin}.
\qed

For $Y\subset X$ as in Proposition \ref{prop_geoext} and a class   $\al$  in $H_n(X,Y)$, consider again  the  support $S:=\spt(\al)$.
For $p\in S$, we  define the space of directions $\Sigma _p S \subset \Sigma _p X$ to be the set
of starting directions of geodesics $py$ contained in $S$.

Due to  Proposition \ref{prop_geoext} any such segment extends within $S$ until $Y$ (in particular, to a uniformly positive length).  By compactness of $\bar S$, this implies that
$\Sigma _p S$ is a compact subset of $\Sigma _p X$. Another application of Proposition \ref{prop_geoext} shows
\[\Sigma_p S=\bigcap_{r>0}\log_p(\dot B_r(p)\cap S)\]
where $\dot B_r(p)$ refers to the punctured ball $B_r(p)\setminus\{p\}$.
We refer to \cite{H} for further properties of $\Sigma _p S$. Here we will only need:

\blem\label{lem_links}
In the notations above, assume that $\Sigma _p S = \mathcal V^+ \cup \mathcal V^-$, with $\mathcal V^{\pm}$ being closed, non-empty  and disjoint.  Then we find $v^+\in \mathcal V^+$ and $v^-\in \mathcal V^-$ with
 \[|v^+,v^-|=\pi.\]
\elem

\begin{proof}
	Assume the contrary and  find a small $\epsilon$, such that $\epsilon<|v^+,v^-|<\pi -\epsilon$
	for all $v^{\pm} \in \mathcal V^{\pm}$.  Choose  arbitrary $v_0^{\pm} \in \mathcal V^{\pm}$
	and $x_0^{\pm} \in S$ lying in the direction of $v_0 ^{\pm}$ from $p$.

	We find a small positive $r>0$, such that for all $q\in  S\cap \partial B_r (p)$ we have
	$| \Sigma _p S, \log _p (q)| \leq  \frac \epsilon 3$. Thus, $S\cap \partial B_r(p)$ is a disjoint union of two compact subsets $K^{\pm}$ such that $|  \mathcal V^{\pm}, \log_p (q)|\leq \frac \epsilon 3 $, for $q\in  K^{\pm}$.

	The geodesics towards $x_0^{\pm}$  provide a contraction of  a neighborhood of  $K ^{\mp}$ in $X$
	   inside  $X\setminus \{p \}$.
	   Thus, a neighborhood of   $S\cap \partial B_r(p) $ is contractible inside $X\setminus \{p \}$, in contradiction to Corollary \ref{cor_fin}.
\end{proof}

\section{Reduction to a cutting Lemma}

Throughout this section
$X$ is a contractible CAT($\kappa$) space
of dimension
two.

\bdfn
Let $\Ga$ be a Jordan curve in $X$ and $[\Ga]$ a generator of  $H_1 (\Gamma)$.
Denote by $\al\in H_2 (X,\Ga)$ the unique element with $\D\al=[\Ga]$.
The support $\spt (\al)  \subset X\setminus \Ga$ will be called \emph{interior} of $\Gamma$.
\edfn
We borrow the term \emph{interior} in this context from \cite{Ricks}.

A closed convex subset $X'$ of $X$ of diameter less than
	$D_{\kappa}$ is contractible. If $\Gamma$ is contained in $X'$,  then, by the definition of support,  the interior of $\Gamma$ is contained in $X'$ and so are all  subsets resulting from the constructions performed in this section. Recall that any Jordan curve $\Gamma$ of length $l < 2D_{\kappa}$ is contained in a closed ball $X'$ of radius $\frac  l 4$, \cite[Prop. 3.20]{Ballmann}.  Since all subsequent considerations concern only such curves, we may always replace $X$ by $X'$ and assume that the diameter of $X$ is at most $\frac l 2 < D_{\kappa}$ to begin with.

 Let $\Ga\subset X$ be a Jordan curve with interior $S$.
We define a {\em cut} of $\Ga$ to be a  geodesic $c$ with $c\subset\bar S$ and $c\cap\Ga=\D c$.
A cut $c$ divides $\Ga$ into two arcs $\Ga^+$ and $\Ga^-$ whose boundaries coincide with $\D c$.
 When performing a cut we obtain two new Jordan curves $\Ga_c^\pm:=\Ga^\pm\cup c$.

 A {\em k-fold (iterated) cut} is defined inductively, where
a 1-fold iterated cut is just a cut and a $k$-fold iterated cut are $2^{k-1}$ cuts performed at each of the Jordan curves resulting from a
$(k-1)$-fold iterated cut.

Iterated cuts  stay in a controlled neighborhood of the original curve:

\blem\label{lem_diamdec}
Let $\Ga\subset X$ be  a
Jordan curve  of diameter at most $\eps <  \frac {D_{\kappa}} 2   $.
Denote by $G_k$
the union of $\Ga$ with all the geodesics from a $k$-fold iterated cut. Then $G_k$ is contained in the $\eps$-neighborhood of $\Ga$.
\elem

\proof
$\Gamma$ is contained in a convex ball  $X'= \bar B_{\epsilon} (x)$, for any $x\in \Gamma$.  Then, by induction,  all iterated cuts of $\Gamma$ are contained in $X'$.
This implies the claim.
\qed

\medskip

For any cut $c$ of a Jordan curve $\Gamma$ of length $< 2D_{\kappa}$ in $X$, both arising Jordan curves $\Gamma _c ^{\pm}$ have length strictly smaller than the length of $\Gamma$.
The same is true for any Jordan curve obtained as a result of an iterated cut. The following cutting lemma is a uniform version of this statement.
We are going to postpone its proof  to the next section.  In this section,  we derive from it the main results of the present paper.

\begin{lem}\label{lem_iteratedcut}
For any  Jordan curve  $\Gamma\subset X$ of length $<2D_{\kappa}$
and   $\epsilon >0$, there exists a $k$-fold cut of $\Gamma$,
 such that all resulting Jordan curves have diameter at most $\eps$.
 \end{lem}

Using this lemma we can now provide:
\proof[Proof of Theorem~\ref{thm_intcurv_intro}]
Thus, let $A\subset X$ be a Lipschitz connected subset of a  contractible CAT($\kappa$)
space $X$
such that  $H_1(A)=0$.

 Since $H_1 (A)=0$,  for any Jordan curve $\Gamma$ contained in $A$, the interior of
 $\Gamma$ is contained in $A$ (by the very definition of support).
 By induction on the number of cuts,
all iterated cuts of $\Gamma$ are contained in $A$.

Denote by $\hat A$ the set $A$ with its intrinsic metric. By assumption $\hat A$ is a intrinsic
	space. Since $A$ is closed,  $\hat A$ is complete \cite[Exercise~1.19]{A_pure}.

	 The identity map $I:\hat A\to A$ is $1$-Lipschitz.  For any Lipschitz curve $\gamma:[a,b]\to A$, the curve $I^{-1} \circ \gamma$ has the same length in $A$ and in $\hat A$. Therefore, for any $1$-Lipschitz map $f:Z\to A$  from an  intrinsic space $Z$, the composition $I^{-1} \circ f$ is $1$-Lipschitz as well. Thus,
	by Proposition \ref{prop_major},  we only need to find for any Jordan curve $\Gamma$ in $A$ of length $<2D_{\kappa}$ a majorization of $\Gamma$
	in $A$ by some CAT($\kappa$) disc $Z$.

	 We fix such a curve $\Gamma$.
 By Lemma~\ref{lem_iteratedcut}, we find an infinite sequence of iterated cuts of $\Ga$ such that the diameters of all occurring Jordan curves
 go  to zero uniformly. Namely, for each  $\eps >0$ there exists $k_\eps\in\N$ such that all the resulting Jordan curves after $k\geq k_\eps$ iterations have diameter at most $\eps$.
 Denote by $G_k$ the union of $\Ga$ with all cuts after $k$ iterations. Then   $(G_k)_{k\in\N}$ forms an increasing  sequence of compact sets in $A$, $G_k\subset G_{k+1}$.
 By Lemma~\ref{lem_diamdec}, we have $G_k\subset N_\eps(G_{k_\eps})$. This implies that the closure $\overline{\bigcup_{k\in\N}G_k}$ is compact.
 Indeed, for all $\eps>0$ we can cover $G_{k_\eps}$ by finitely many $\eps$-balls. Hence $\overline{\bigcup_{k\in\N}G_k}\subset \bar N_{\eps}(G_{k_\eps})$
 is covered by finitely many $2\eps$-balls.

Since all Jordan curves which arose from our cutting process have diameter at most $\eps$, Lemma~\ref{lem_controlledmaj} ensures that they can be majorized within their closed
$\eps$-neighborhoods. Using Lemma~\ref{lem_glumaj}, we can then inductively glue these majorizations to obtain a majorization $f_k:Z_k\to X$ of $\Ga$ with image
contained in $\bar N_{\eps}(G_{k})$.
Precomposing with a majorization of $\partial Z_k \subset Z_k$, we may assume that $Z_k$ is a convex region in $M^2_\kappa$.
The sequence $(Z_k)_{k\in\N}$ subconverges with respect to Hausdorff distance to a convex region $Z_\infty$.
 Thus we obtain a partial limit $f_\infty:Z_\infty\to X$ which is 1-Lipschitz and has image in $\overline{\bigcup_{k\in\N}G_k}\subset A$.
 To see that it is a majorization of $\Ga$, we note
 \[\mathcal \ell(\Ga)\leq \mathcal \ell(\D Z_\infty)\leq\liminf \mathcal \ell (\D Z_k)=\mathcal \ell (\Ga).\]
 Hence the boundary $\D Z_\infty$ is mapped in an arc length preserving way by $f_\infty$.
\qed

\medskip

From the  above proof we now  derive
\proof[Proof of Corollary~\ref{cor_local}]
As before, denote by $\hat A$ the set $A$ with the induced intrinsic metric.
Let  $x\in A$ be an arbitrary point.  Since $Y$ has curvature bounded above by $ \kappa$,
we find  a small  ball $\bar B_r(x)$ around $x$ in $Y$ which is CAT$(\kappa)$ and such that $r < \frac {D_{\kappa}} 2$.

Since $A$ is locally simply connected at $x$,  we find some $s<r $, such that the inclusion $\bar B_s(x) \cap A\to B_r(x)\cap A$ induces a trivial map on $\pi _1$, hence also on  $H_1$.

Denote by $P$ the set  of all points in $\bar B_s(x) \cap A$ which are connected
to $x$ by a Lipschitz curve inside $\bar B_s (x) \cap A$. Let $\hat P$ be the set $P$ equipped with the induced intrinsic metric.    Note that the $\frac s 4$-ball in $\hat P$ around $x$ coincides with the $\frac s 4$-ball around $x$ in $\hat A$. Thus, it suffices to prove that $\hat P$ is CAT$(\kappa)$.
The space $\hat P$ is a complete intrinsic space, by the same argument as in the solution of \cite[Exercise~1.19]{A_pure}.
Assume that $\Gamma$ is a Jordan curve in $\hat P$ of length $<2D_{\kappa}$.  Then $\Gamma$ defines a trivial element in $H_1 (A\cap B_r(x))$, thus, the interior of $\Gamma$ within $\bar B_r(x)$  is contained in $A\cap \bar B_r(x)$.    So, any cut of $\Gamma$ (considered within the CAT($\kappa$)-space $\bar B_r(x)$) is
contained in $A\cap \bar B_r(x)$. Since a cut is a geodesic, it is also contained in the convex ball $\bar B_s(x)$.

 Repeating the argument,  any iterated cut of $\Gamma$ is contained in $A\cap \bar B_s(x)$.  As in the proof of Theorem~\ref{thm_intcurv_intro}, we  now find  a majorization of $\Gamma$ within $A\cap \bar B_s(x)$
by a  convex subset of $M^2_{\kappa}$.

The image of this majorization lies in $P$ and provides a  majorization of $\Gamma$ within $\hat P$. By  Proposition \ref{prop_major}, this implies that $\hat P$ is CAT$(\kappa)$ and finishes the proof.
\qed

\medskip

In the case $\kappa=0$, the theorem of Cartan--Hadamard  yields
  Corollary~\ref{cor_asph_intro}.

\section{Proof of the cutting Lemma}

\subsection{Reduction to  existence of  essential cuts}
   For $\delta >0$,  a $k$-fold cut of $\Gamma$ will be  called  {\em $\delta$-essential}, if all $2^k$ resulting  Jordan curves $\Ga_i$  satisfy
$ \ell (\Ga_i)\leq(1-\delta)\cdot \ell (\Ga)$.

Lemma~\ref{lem_iteratedcut} is a direct consequence of the following:

\begin{prop}\label{prop_esscut}
For every $\eps _0>0$,  $\kappa\in\R$ and $2D_{\kappa} >l_0 >0$ there exists a positive constant $\delta=\delta (\eps_0, \kappa, l_0)$ such that the following holds:

	Let $\Ga$ be a Jordan curve of length $l \leq l_0$ in a CAT($\kappa$)
	space $X$.
	 If the diameter of $\Gamma$ is at least
	 $\eps_0$, then $\Ga$
	admits a $\delta$-essential
	 2-fold cut.
\end{prop}

Assuming Proposition  \ref{prop_esscut}, we now provide
\proof[Proof of Lemma~\ref{lem_iteratedcut}]
Let $l$ be the length of $\Ga$. Choose $\delta= \delta(\eps, \kappa , l)$  in Proposition~\ref{prop_esscut}.
If $\Gamma$ has diameter $\leq \epsilon$ there is nothing to prove.
Otherwise, we apply Proposition~\ref{prop_esscut} and obtain a $\delta$-essential
2-fold cut of $\Gamma$ resulting in four new Jordan curves.

We proceed inductively as follows. At each step we consider the Jordan curves produced by the previous step and split them into two groups depending
on their diameter. If a Jordan curve has diameter  at most $\epsilon$, we perform an arbitrary 2-fold cut, thereby keeping the diameter bound.
If, on the other hand,
a Jordan curve has diameter larger than $\epsilon$, then, since its length is less than $l$, we can apply Proposition~\ref{prop_esscut} to make a
$\delta$-essential 2-fold cut.

Thus, after $k$ steps, the diameters  of the arising Jordan curves are at most
$ \max \{\epsilon, (1-\delta )^k \cdot l  \}$. For sufficiently large $k$,
this number is at most $\epsilon$ as required.
\qed

\subsection{Setting for finding essential cuts}
In order to prove  Proposition~\ref{prop_esscut}, we fix
	  $\kappa$,  $l_0 < 2D_{\kappa}$   and $\epsilon _0 $.
	By scaling, we may assume $\kappa=1$.  Making $\epsilon_0$ smaller, we  may   assume $\epsilon  _0  <  \frac 1 8$.

	 Henceforth we fix a CAT(1) space $X$ and a Jordan curve $\Gamma$ in $X$ of length  $l \in (\epsilon _0 , l_0]$.  The curve $\Gamma$ is contained in a ball of radius at most $\frac l 4$ and replacing $X$ by this ball, we may assume that $X$ has diameter at most $\frac l 2$.
	     Let $S$ denote  the interior of $\Gamma$.
	We are going to construct a 2-fold $\delta$-essential cut of $\Gamma$ for some
	 $\delta$ depending only on
	 $l_0$ and $\epsilon _0 $.

The proof will be divided into two cases, depending on
the existence of
 long cuts:
 For $\epsilon >0$, we  say that  $\Gamma$ is \emph{$\epsilon$-degenerated}
if   all cuts of $\Gamma$ are shorter than $2\epsilon \cdot l$.

\subsection{The non-degenerated case}
In this case we proceed as follows.  First we are going to verify the following simple claim about spherical triangles:

There exists  a positive constant $\rho = \rho (l_0)$  such that
for any spherical triangle $\triangle(x,y,z)\subset\mathbb S^2$  with $|y,z| \leq\pi /2$,  $|x,y|\leq    \frac {l_0 }  2$, $\angle_y(x,z) \geq \frac \pi 2$ we have
\[|x,y| \leq |x,z|+ |z,y| - \rho \cdot | z,y| .\]

If $|x,y| \leq \frac \pi 2 $, we can take $\rho =1$.

For $\frac {l_0} 2\geq |x,y| \geq \frac \pi 2$,  we may assume  $\angle_y(x,z)=\frac \pi 2$ and $|x,y|=\frac {l_0} 2$.
The claim follows by compactness and  the fact that for sufficiently small $|x,y|$ the statement is true  by
the first variation formula.

Using this $\rho =\rho (l_0)$  we can now state:

\blem\label{lem_crosscut}
If, for some $\eps\in(0,\frac 1 8)$, the curve   $\Gamma$ is not  $\epsilon$-degenerated then there is  a $(\rho\epsilon)$-essential 2-fold cut of $\Gamma$.
\elem

\proof
Let $c$ be any cut of length at least $2\epsilon \cdot l$, which exists by assumption.
Denote by $z_1, z_2 \in \Gamma$ the endpoints of $c$.
Denote by $\Gamma ^+ _c$ and $\Gamma ^- _c$ the arising Jordan curves with interiors $S^{\pm}$.
 Let $m$ denote the midpoint of $c$.

In $S^+$ we consider a sequence of points $x_k$ converging to $m\in c \subset \Gamma ^+$, which is possible by Lemma~\ref{lem_acc}.
Denote by $\bar x_k\in\Ga^+$ the closest point to $x_k$. Then $\bar x_k$ converges to $m$, and  taking $k$ large enough, we may assume
that $\bar x_k \in c$.
Using the geodesic extension property, Proposition~\ref{prop_geoext}, we can extend the segments $\bar x_k x_k$  to a segment $c_k= \bar x_k y_k$  such that $x_k y_k$ is contained in the closure $\bar S^+$. Since $\bar S ^+$ is compact, we can pass to a partial limit $c'$ of $c_k$, which is a geodesic segment $my$ starting in $m$ and contained in $\bar S^+$.
By the upper semi-continuity of angles,  both angles  enclosed between  $c$ and $c'$ are at least $\frac \pi 2$, since the same is true for $c_k$ and $c$ at
$\bar x_k$.   Restricting
$c'$ to the first intersection point $y^+ \neq m$ of $my$ with $\Gamma ^+$ we obtain a cut
$c^+ =my^+$ of $\Gamma ^+$ starting at $m$ and enclosing  angles at least $\frac \pi 2 $ with $c$.

In the same way, we construct a cut $c^-=my^-$ of $\Gamma ^-$  starting at $m$ and
enclosing angles at least $\frac \pi 2$ with $c$.
The geodesics $c, c^+,c^-$ provide a 2-fold cut of $\Gamma$.  We claim that this 2-fold cut is $(\rho \epsilon)$-essential.

\includegraphics[scale=0.6,trim={0cm 2cm 0cm 0cm},clip]{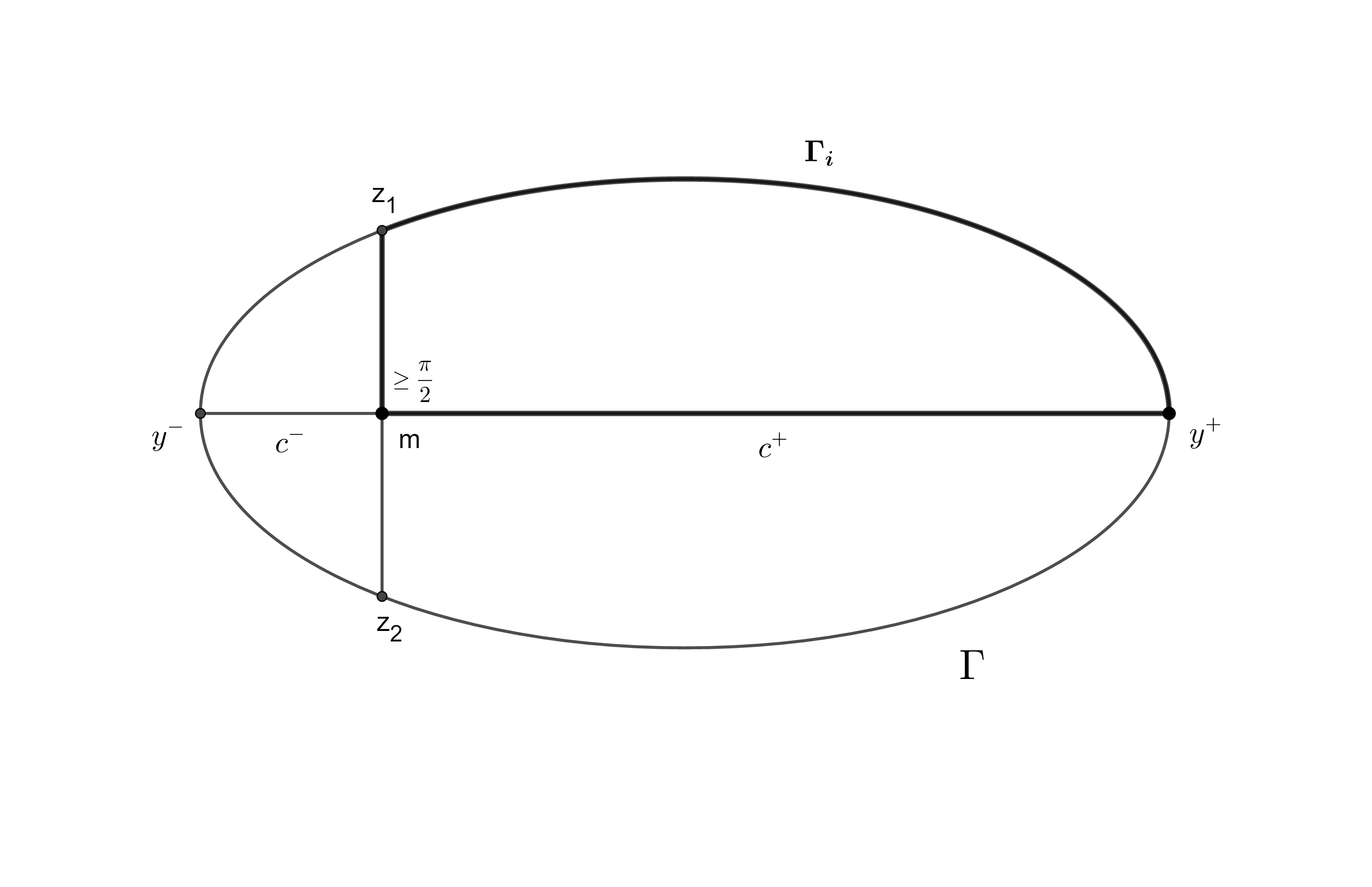}

In order to see this, consider one of the four Jordan curves $\Gamma _i$. Without loss of generality, we may assume that $\Gamma _i$ consists of the geodesics $z_1m$, $my^+$ and the part of $\Gamma$ between $y^+$ and $z_1$.  Thus, $\Gamma _i$ arose  from the cut $c^+$ of $\Gamma ^+$. The difference of lengths
$$\ell (\Gamma ^+) -\ell (\Gamma _i)$$
is at least as large as the triangular defect
$$|y^+,z_2|+|z_2,m| -|m,y^+|$$
which by construction of $\rho$ is at least  $$\rho \cdot  |z_2, m|  \geq \rho \cdot \epsilon \cdot l \;.$$
Therefore, we deduce
$$\ell   (\Gamma) -\ell (\Gamma _i) \geq \ell (\Gamma ^+) -\ell (\Gamma _i) \geq \rho \cdot \epsilon \cdot l \;.$$
This finishes the proof.
\qed

\subsection{The degenerated case} This case is technically more complicated and requires several steps.  In the rest of this subsection, we fix
$\delta:= \frac {\epsilon _0} {1000 \cdot \pi}$  and   assume that
$\Gamma$ is $\delta$-degenerated.   Since $l <2\pi$, we have $5\cdot\delta \cdot l < \frac {\epsilon_0} {100}$.

 We note that the closure  $\bar S$ of $S$ is Lipschitz connected: $\Ga$ has finite length and any point in $S$ is connected to some  point  on $\Gamma$ by an (extrinsic) geodesic completely contained in $\bar S$, due to Proposition \ref{prop_geoext}.

 Denote by $\hat S$ the set $\bar S$ equipped with the induced intrinsic metric.
If we consider $\Gamma$ as a subset of $\hat S$ we will denote it by $\hat \Ga$.  The identification $\hat \Gamma \to \Gamma$ is a length preserving, $1$-Lipschitz homeomorphism.

 We observe:
$\hat \Gamma$ is $ (\delta \cdot   l)$-dense in $\hat S$.

Otherwise, we could find a point $x\in \hat S$ at distance at least $ \delta \cdot l $ from $\hat\Ga$.  By Proposition~\ref{prop_geoext}, we find a cut
$c$ of $\Gamma$ through the point $x$. By assumption on $x$, the cut $c$ has length at least $2\delta \cdot l$, in contradiction to our degeneracy assumption.

Using this observation, we are going to control the \emph{absolute filling radius} of $\hat \Gamma$.   Recall
that the absolute filling radius of $\hat \Gamma$ is
 the greatest lower bound of numbers $r$ such that $\hat \Gamma$ embeds isometrically as an $r$-dense subset into a metric space $Y$,
$\iota: \hat \Gamma \hookrightarrow Y$, and $\iota_\ast [\hat \Ga] =0\in H_1(Y)$.

We will use the following basic observations about absolute filling radii:
\begin{itemize}
\item The absolute filling radius of  $\mathbb S^1$ with its intrinsic metric  is $\frac{\pi}{3}$, \cite{Ka_filling}.
\item If $f:\hat \Gamma \to \tilde \Gamma$ is an $L$-Lipschitz map of degree one, then the absolute filling radius of $\tilde \Ga$
is at most $L$ times  the filling radius of $\hat \Ga$, \cite[p.~8]{G_filling}.
\item If $|\cdot , \cdot|_k$  is a sequence of metrics on $\hat \Gamma$ converging uniformly to $|\cdot, \cdot|_{\hat \Gamma}$
  then the absolute filling radii of $(\hat \Gamma, |\cdot , \cdot|_k)$ converge to the absolute filling radius of $(\hat \Gamma ,|\cdot, \cdot|_{\hat \Gamma})$, \cite[Proposition 9.34]{LMO}.
\end{itemize}

Under our degeneracy   assumption we can now show:

\blem\label{lem_fillradbound}
The Jordan curve $\hat \Gamma \subset \hat S$ has absolute filling radius at most $ \delta  \cdot l$.
\elem

\proof

For every $k\in\N$, consider the neighborhood $W_k =  N_{\frac{1}{k}}(\bar S)$ of $\bar S$.  Due to  Lemma~\ref{lem_acc},   the curve $\Gamma$ bounds a chain in $W_k$.  This relative cycle can be represented  by a continuous map $f:\Sigma  \to W_k$, where $\Sigma$
is a smooth Riemannian  surface  with one boundary curve which is mapped  by $f$ onto $\Gamma$ in a length preserving way, \cite[p.~109]{Hatcher}.  Due to the straightening of simplices  mentioned above, \cite[Section~6.1]{KleinerLeeb}, we may assume that $f$ is Lipschitz continuous after perturbation.

The set $W_k$ is also Lipschitz connected. We denote by $\hat W_k$ the set $W_k$ with its intrinsic metric and by $\Gamma _k$, the curve $\Gamma$ as a subset of $\hat W_k$.

Since $f$ is Lipschitz continuous, it remains  continuous as a map $f:\Sigma \to \hat W_k$.
Thus, $\Gamma _k$ induces a trivial $1$-cycle in $H_1 (\hat W_k)$. Since any point of $\bar S$ is connected to $\Gamma$ with a curve of length at most $\delta \cdot l$,
the curve $\Gamma _k$ is $(\delta\cdot l + \frac 1 k)$-dense   in $\hat W_k$.  Therefore, the filling radius of $\Gamma _k$ is at most $(\delta\cdot l + \frac 1 k)$.

In order to reach the conclusion, we only need to show that $| \cdot ,\cdot |_{\Gamma_k}$ converge
uniformly to $|\cdot ,\cdot |_{\hat \Gamma}$. By  the uniform compactness of $\Gamma _k$ it suffices to prove that for any $p, q\in \Gamma$ the distances
$|p,q|_{\Gamma _k}= |p,q|_{W_k}$ converge to $|p,q|_{\hat \Gamma}  =|p,q|_{\hat S}$.

Clearly, $|p,q|_{\hat W_k} \leq  |p,q| _{\hat S}$.

 On the other hand,   consider a sequence of curves $\eta _k$ in $W_k$ connecting $p,q$, parametrized by arclength and such that $\lim \limits_{k\to\infty} \ell(\eta _k) = \lim\limits_{k\to\infty} |p,q| _{\hat W_k}$.  By compactness of $\bar S$, we find a subsequence of $\eta_k$ which  converges pointwise to a curve $\eta$ in $ \bar S$.
 Then $\ell (\eta) \leq  \lim \limits_{k\to\infty} \ell(\eta _k) $, hence
\[\lim\limits_{k\to\infty} |p,q|_{\hat W_k}    \geq |p,q|_{\hat S} \,.\]
This  finishes the proof of the lemma.
\qed

\medskip

We are going to show  that the $\delta$-degenerated $\hat\Ga$  ``comes close to itself''.

	Consider points $p, q \in \hat \Ga$ realizing the diameter $d$ of $\hat \Gamma$.
	Since the diameter of $\Gamma$ does not exceed the diameter of $\hat \Gamma$,  by our assumption $d >\epsilon_0$.

	Denote by $\Gamma ^{\pm}$ the arcs of $\hat \Gamma$ defined by the points $p$ and $q$.
Denote by $\Gamma ^{\pm} _0$ the set of all points in $\Gamma ^{\pm} $ that have $\hat S$-distance at least $\frac d {3}$ to $p$ and $q$.     We claim

\begin{lem}\label{lem_degone}
There exist points $x^+\in \Gamma^+ _0$ and $x^- \in \Gamma ^- _0$  with $|x^+,x^-| _{\hat S} \leq 4\cdot \delta \cdot l$.
\end{lem}

\proof
Assume the contrary.
	The function $\hat f(x):=|p,x|_{\hat S}$ defines $1$-Lipschitz maps of both arcs $\Gamma ^{\pm}$
	onto the interval $[0,d]$.  Moreover, $\hat f$ sends the points $p,q$ onto the ends of the interval and is a degree one  map  of $\Gamma ^{\pm}$ onto $[0,d]$ modulo endpoints.

	Hence also the maps  $\tilde f:\Gamma ^{\pm} \to [\frac d {3} , \frac {2d} {3}],$ defined as composition of $\hat f$ and the closest-point projection from $[0,d]$ to
	$[\frac d {3} , \frac {2d} {3}]$, are $1$-Lipschitz and have degree one.

Consider the constant speed parametrizations $\eta ^{\pm}  : [\frac d {3} , \frac {2d} {3}] \to \mathbb S^1$  of the upper and lower hemi-circle, respectively. Let
$ f:\hat \Gamma \to \mathbb S^1$ be  defined on $\Gamma ^{\pm}$  as $\eta ^{\pm}\circ \tilde f$.  By construction, the map $ f$ has degree one  and its restrictions to $\Gamma^+$ and to $\Gamma ^-$ are both $(\frac {3}  {d} \cdot \pi)$-Lipschitz continuous.

For all  $x^+ \in \Gamma ^+$ and $x^- \in \Gamma ^-$, either one of the points is not in $\Gamma ^{\pm} _0$ and then
$$\frac {|f(x^+) , f(x^-)|} {|x^+,x^-| _{\hat S}}  \leq \frac {3}  {d} \cdot \pi \,,$$
by the $1$-Lipschitz property of $\tilde f$. Or, otherwise,  the distance between $x^+$ and $x^-$ is at least $4\cdot \delta \cdot l$ and then
$$\frac {|f(x^+) , f(x^-)|} {|x^+,x^-| _{\hat S} }  \leq \frac  {\pi} {4\cdot \delta \cdot l} \,.$$
By assumption  on $\delta $
 we have  $\frac 3 {d} \leq \frac 1 {4\cdot \delta \cdot  l}$.
 Thus, the map $f$
 is $ \frac \pi {4\cdot \delta \cdot  l}$-Lipschitz.

Due to   Lemma~\ref{lem_fillradbound} and the properties of the filling radius listed above,  we deduce $$\frac \pi {4\cdot \delta \cdot  l}  \cdot \delta \cdot l     \geq   \frac \pi 3 \;.$$
This is a contradiction, finishing the proof.
\qed

\medskip

Now we will find a ``good'' cut near the points provided by the above lemma:

\bcor\label{cor_degone}
The curve $\Gamma$ admits a $\delta$-essential cut.
\ecor

\proof
We continue to use notations introduced prior to Lemma~\ref{lem_degone}.
Consider points $x^{\pm}  \in \Gamma ^{\pm} _0$ provided by  Lemma~\ref{lem_degone}.
We find a curve   $ \hat \eta$ in $ S$ connecting $x^+$ with $x^-$ of length  $<5\cdot \delta \cdot l$.
 Let $x^+_0$ be the last intersection point of $\hat \eta$ with $\Gamma ^+$ and let $x^- _0$ be the first intersection point of $\hat \eta$ with $\Gamma ^-$.  Denote by
$\eta$  the part of $\hat \eta$ between $x^+_0$ and $x^- _0$.
\vspace{0.5cm}
\begin{center}
\includegraphics[scale=0.28,trim={0cm 0cm 0cm 0cm},clip]{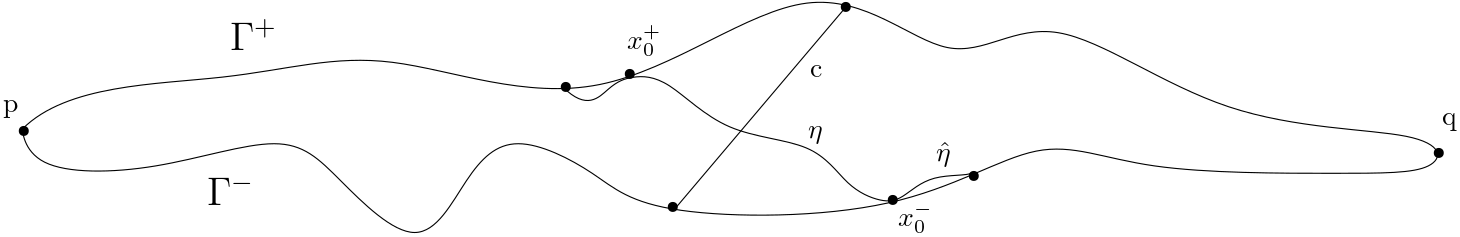}
\end{center}
\vspace{1cm}

Thus,  the length of   $\eta$ is $<5 \cdot \delta \cdot l$ and, by construction,
 the distances from $x^{\pm} _0$ to $p$ and $q$  are at least
$$\frac d {3}- 5\cdot \delta  \cdot l  \geq \frac d {3} - \frac d { 100} \geq \frac d 4 \,.$$

Consider the set $\mathcal C$ of all cuts $c$ of $\Gamma$ which contain a point on $\eta$. Let $\bar {\mathcal C}$ be the set of all geodesics which can be obtained as a limit of a sequence in $\mathcal C$.

Any geodesic $\hat c \in \bar {\mathcal C}$ connects two points on $\Gamma$ and intersects  $\eta$. By the degeneracy assumption, $\hat c$ has length at most $2\cdot \delta \cdot  l$.
Thus, any point on $\hat c$ has distance  from $p$ and $q$ at least equal to
$$\frac d {4} - 7\cdot \delta \cdot  l   \geq \frac d {5} \,.$$

Assume that there exists an element $\hat c \in \bar {\mathcal C}$ which contains points from $\Gamma ^+$ and $\Gamma ^-$ simultaneously.  Then
 a subsegment $c$ of $\hat c$
 		 is a cut of $\Gamma$ and has one endpoint on $\Gamma ^+$ and the other endpoint on $\Gamma^-$.  Then $c$ subdivides $\Gamma$ into two arcs, one of which contains $p$ and the other contains $q$, hence both of them are at least   $\frac {2d} 5$ long.

Thus, both Jordan curves $\Ga_c^\pm$ arising from the cut $c$ have lengths at most
$$l  -  \frac {2d} {5} + 2\cdot \delta \cdot l \leq l-\delta \cdot l\;.$$
Therefore,
$c$ is $\delta$-essential and the proof would be complete.

Therefore,  assuming that the corollary does not hold,  we infer that no  $\hat c\in \bar {\mathcal C}$  simultaneously  intersects $\Gamma ^+$ and $\Gamma ^-$.   Denote by $\bar {\mathcal C}^{\pm}$ the subsets of elements of $\bar {\mathcal C}$ which contain points in $\Gamma ^+$, respectively in $\Gamma ^{-}$.  As we have observed, our assumption implies  that  the sets $\bar {\mathcal C}^{\pm}$ are disjoint.  By definition, both these sets are closed under convergence, and
we have $\bar {\mathcal C}  = \bar {\mathcal C} ^+ \cup \bar {\mathcal C} ^-$.

By Proposition~\ref{prop_geoext}, there exists an element $c \in \mathcal C$ through each point in $\eta \setminus \{ x_0 ^{\pm}  \} $.   Denote by $K^{\pm}$  the set of points on $\eta \setminus \{ x_0 ^{\pm}  \} $, which lie on some segment in  $\bar {\mathcal C} ^+$ and $\bar {\mathcal C} ^-$, respectively.   Hence, $\eta \setminus \{ x_0 ^{\pm}  \}  =K^+\cup K^-$.
  Since $\bar {\mathcal C }^{\pm}$ is closed, the sets $K^{\pm}$ are closed in the connected set $\eta \setminus  \{ x_0 ^{\pm}  \} $.

We claim that the sets $K^{\pm}$ are non-empty. Indeed, consider points $y_m  \neq x_0 ^+$ on $\eta$  converging to $x_0 ^+$. Choose cuts $c_m $ of $\Gamma$ through $y_m$.   If $c_m \in \bar {\mathcal { C}}^-$ for all $m$, then any  partial limit $\hat c$ of $c_m$  is contained in
	$\bar {\mathcal { C}}^-$. But $\hat c$ contains $x_0^+$, in contradiction  to the disjointness of $\hat c \in \bar {\mathcal C} ^-$ and  $\Gamma  ^+$.   Thus, $c_m \in  \bar {\mathcal  C}^+$, for large $m$. Therefore, $K^+ $ is not empty. Similarly, $K^-$ is non-empty as well.

The connectedness of $\eta \setminus  \{ x_0 ^{\pm}  \} $  implies that $K^+\cap K^-$ is not empty.
Consider an arbitrary point  $z \in  K^+\cap K^-$.

Denote by $\mathcal V^{\pm}$ the set of  all vectors $v\in \Si _z X$ tangent to  an element in
 $\bar {\mathcal C} ^{\pm}$.   By construction,
  $\mathcal V^{\pm}$ are both non-empty and their union  $\mathcal V^+ \cup \mathcal V^-$ is the space of directions $\Sigma _z S$.
  Since the sets $\bar {\mathcal C} ^{\pm}$
 are closed, the sets $\mathcal V^{\pm}$ are closed in $\Sigma _z X$.  Since  $\bar {\mathcal C} ^{+} \cap\bar {\mathcal C} ^{-}$ is empty,    the intersection of $\mathcal V^+$ and $\mathcal V^-$ is empty and
 no direction
 $v\in \mathcal V^+$ has an antipodal direction in $\mathcal V^-$.
 This contradicts Lemma~\ref{lem_links} and finishes the proof.
\qed

\subsection{Conclusion}
Now we finish the proof of  the existence of essential  iterated cuts.

\proof[Proof of Proposition~\ref{prop_esscut}]
Set $\delta = \frac {\epsilon_0} {1000 \cdot \pi}$. If $\Ga$ is $\delta$-degenerated  we obtain a $\delta$-essential cut from
Corollary~\ref{cor_degone}. If $\Gamma$ is not $\delta$-degenerated, we apply Lemma~\ref{lem_crosscut} and    obtain a
 $(\rho \cdot \delta)$-essential 2-fold cut, where $\rho$
depends only on $l_0$.
\qed

\section{Lipschitz homotopy groups}

In this final section we provide proofs for the  applications of the main result.

\proof[Proof of Theorem~\ref{thm_lipasph_intro}]
Let $A$ be a subset of a two-dimensional non-positively curved space $Y$, let $n \geq2$ be fixed and let $f:\mathbb{S}^n \to A$ be a Lipschitz map.
In order to prove that $f$ is contractible in $A$, we may replace $A$ by the compact image $K=f(\mathbb{S}^n)$ and assume that $A$ is compact.
Rescaling the metric on $Y$ by a factor, we may assume that $f$ is $1$-Lipschitz if $\mathbb{S}^n$ is considered with respect to its intrinsic metric.

Using that  $\mathbb{S}^n$ is simply connected we can lift $f$ to a map  $\tilde f:\mathbb{S}^n \to X$ into the universal covering $X$ of $Y$.  Note that $\tilde f$ is still $1$-Lipschitz  and $X$ is  CAT(0),
 by the theorem of Cartan--Hadamard.   Once we can contract $\tilde f$ in its image, we can also contract $f$.  Thus we may assume that $Y=X$ is CAT(0), that $A$ is compact and that $f$ is $1$-Lipschitz with respect to the intrinsic metric on $\mathbb{S}^n$.

 For $\eps>0$ we can cover the compact set $A$ by finitely many closed $\eps$-balls.
Denote by $A_\eps$ the union of these balls. Then $A_\eps$ is a closed and Lipschitz connected subset of $X$.
Moreover, $A_\eps$ is locally contractible since the union of those balls which contain a fixed point is star shaped. Denote by $\hat A_\eps$
the set $A_{\eps}$ with its the intrinsic metric.
It follows from Corollary~\ref{cor_local} that $\hat A_\eps$ is non-positively curved. Note that $f$ is still 1-Lipschitz as a map from the intrinsic space $\mathbb{S}^n$ to $\hat A_\eps$.
Denote by $\pi_\eps:Y_\eps\to\hat A_\eps$ the universal cover. We lift $f$ to a 1-Lipschitz map $\tilde f:\mathbb{S}^n\to Y_\eps$.
Then $\tilde f$ is $\frac{\pi}{2}$-Lipschitz when considering $\mathbb{S}^n$ as a subset of $\R^{n+1}$.

Since $Y_\eps$ is CAT(0), Kirszbraun's theorem
\cite{LaSch_kirsz,AKP_kirsz} implies that $\tilde f$  has a $\frac{\pi}{2}$-Lipschitz extension
$\tilde F:B^{n+1}\to Y_\eps$ where $B^{n+1}$ denotes the closed unit ball in $\R^{n+1}$.
Projecting to $A_{\eps}$, we extend $f$  to a $\frac{\pi}{2}$-Lipschitz map
 $F:=\pi_\eps\circ\tilde F :B^{n+1} \to \hat A_{\eps}$.

 Hence $F$ is also a $\frac{\pi}{2}$-Lipschitz extension of $f$ as a map $B^{n+1} \to A_{\eps}$,
 where $A_{\eps}$ is equipped with the induced metric from $Y$ (and not with its intrinsic metric).

Now we choose a sequence $\eps_k\to 0$. We obtain a nested sequence $A_{\eps_k}$ of closed sets with $A=\bigcap_{k\in \N} A_{\eps_k}$ and
$\frac{\pi}{2}$-Lipschitz maps $F_k:B^{n+1}\to\hat A_{\eps_k}$ filling $f$. All the maps $F_k$ are $\frac{\pi}{2}$-Lipschitz as maps to $Y$.
Since $K$ is compact, we obtain a partial  limit
$F_\infty:B^{n+1}\to Y$ which is a $\frac{\pi}{2}$-Lipschitz map  extending  $f$ and has image in $A$.
\qed

\proof[Proof of Corollary~\ref{cor_opasph_intro}.]
Let $Y$ be a two-dimensional space of non-positive curvature, let $A\subset Y$ be a neighborhood retract and let $f:\mathbb{S}^n \to A$
be continuous for some $n\geq 2$.  By assumption, we find a retraction $r:U\to A$ of an  open neighborhood $U$ of $A$ and it suffices to find a continuous  filling  $F:B^{n+1} \to U$ of $f$.

 Thus, we may assume that $A=U$ is open in $Y$. By staightening simplices, we see that $f$ is homotopic to a Lipschitz map $\hat f:\mathbb{S}^n \to U$, \cite[Section~6.1]{KleinerLeeb}.   By Theorem~\ref{thm_lipasph_intro}  the map $\hat f$ is contractible in its image, hence $f$ is contractible in $U$.
\qed

\medskip

In the proof of Proposition~\ref{prop_coarsetop} we will use the concept of ultralimits of metric spaces with respect to a non-principal ultrafilter $\om$.
For precise definitions and properties we refer the reader to \cite{KleinerLeeb}, \cite[Chapter~1.5]{BH},  \cite[Chapter~10]{DK_ggt}.

\proof[Proof of Proposition~\ref{prop_coarsetop}]

Suppose for contradiction that the claim is wrong for some $n$ and $\eps$.
Then we find a sequence $L_k\to\infty$ and a sequence of $L_k$-Lipschitz maps $f_k:\mathbb S^{n}\to X$ with images in $A$ which do not bound Lipschitz balls
in $N_{\eps L_k}(A)$.

 Let $p$ be a base point in $\mathbb S^{n}$ and set $x_k:=f_k(p)$. Then we rescale $X$ by $\frac{1}{L_k}$ and denote the resulting space by $X_k$.
Note that $f_k$ is 1-Lipschitz as a map to $X_k$. We pass to ultralimits $(X_\om,x_\om)=\wlim(X_k,p_k)$ and $f_\om:\mathbb S^{n}\to X_\om$. Then $f_\om$ is 1-Lipschitz and by assumption, $X_\om$ is a CAT(0) space of dimension at most two. By the proof of Theorem~\ref{thm_lipasph_intro},
$f_\om$ bounds a $\frac{\pi}{2}$-Lipschitz ball $F_\om$ in its image.
Choose a finite  $\frac{\eps}{2\pi}$-dense set $T$ in the open unit ball $B^{n+1}$ in $\R^{n+1}$.
Define  an extension $\tilde F_k:\mathbb S^{n}\cup T\to X_k$ of $f_k$ such that $\wlim \tilde F_k=F_\om|_{S^n\cup T}$.
Then for $k$ large enough, all $\tilde F_k$ are $\pi$-Lipschitz (where $\mathbb S^{n}\cup T$ carries the induced metric from $\R^{n+1}$). By Kirszbraun's theorem~\cite{LaSch_kirsz,AKP_kirsz}, we can extend $\tilde F_k$ to a $\pi$-Lipschitz map $F_k:\bar B^{n+1}\to X_k$.
Again, for $k$ large enough, $F_k(T)$ lies at distance $<\frac{\eps}{2}$ from $A$ and therefore the image of $F_k$ lies at distance $<\eps$ from $A$. Contradiction.
\qed

\bibliographystyle{alpha}
\bibliography{2D}

\Addresses

\end{document}